\documentclass{amsart}

\usepackage{amsmath}
\usepackage{amssymb}
\usepackage{amsthm}

\newcommand{\Z}{{\bf Z}}

\newcommand{\F}{{\bf F}}
\newcommand{\R}{{\bf R}}

\newtheorem{theorem}{Theorem}

\newtheorem{question}{Question}
\newtheorem{prop}{Proposition}[section]

\begin{document}

\title[Abelian varieties with a specified characteristic
polynomial]{Abelian varieties over finite fields with a specified
characteristic polynomial modulo $\ell$}

\author{Joshua Holden}

\address{Department of Mathematics, Duke University, Box 90320\\
Durham, NC 27708-0320, USA}

\email{holden@math.duke.edu}

\begin{abstract}
We estimate the fraction of isogeny classes of abelian varieties over
a finite field which have a given characteristic polynomial $P(T)$
modulo $\ell$.  As an application we find the proportion of isogeny
classes of abelian varieties with a rational point of order $\ell$.
\end{abstract}

\maketitle

\section{Introduction}

Let $\F$ be a finite field of characteristic $p$ and order $q$, and
$\ell$ a prime not equal to $p$.
Let
\begin{equation*}
P(T)=(T^{2g}+q^g)+a_1(T^{2g-1}+q^{g-1}T)+ \cdots + a_{g-1}(T^{g+1}+q
T^{g-1})+a_g T^g
\end{equation*}
be a polynomial.  The goal of this paper is to estimate the number of
isogeny classes of abelian varieties over $\F$ of dimension $g$ for
which the characteristic polynomial for the action of Frobenius is
congruent to $P(T)$ modulo $\ell$.

The initial motivation for this problem came from the following
question, posed in~\cite{Holden99} and related to the Fontaine-Mazur
Conjecture for number fields:

\begin{question} \label{falseconj1}
Let $k$ be a function field over a finite field $\F$ of characteristic
$p$ and order $q$, and $\ell$ a prime not equal to $p$.  Let $K =
k\F_{\ell^\infty}$ be obtained from $k$ by taking the maximal
$\ell$-extension of the constant field.  If $M$ is an unramified
$\ell$-adic analytic $\ell$-extension of $k$, and $M$ does not contain
$K$, must $M$ be a finite extension of $k$?
\end{question}

In general the answer to Question~\ref{falseconj1} is no, with
examples due to Ihara (\cite{Ihara83}) and to Frey, Kani, and
V\"olklein (\cite{FKV}).  However, the following theorems were proved
in~\cite{Holden99}:

\begin{theorem}[Theorem~2 of~\cite{Holden99}] \label{thm2}
Let $k_0$ be a function field over a finite field of characteristic
$p$, and let $k$ be a constant field extension.  Let $\ell$ be a prime
not equal to $p$.  If $\ell$ does not divide the class number $P(1)$ of
$k_0$, then any everywhere unramified powerful (\emph{a fortiori}
uniform) pro-$\ell$ extension of $k$, Galois over $k_0$, with no
constant field extension, is finite.
\end{theorem}

\begin{theorem}[See Corollary~4.11 of~\cite{Holden99}] \label{symm-free}
Let $k_0$ be a function field over a finite field of characteristic
$p$, and let $k$ be a constant field extension.  Let $\ell$ be a prime
not equal to $p$.  Let $P(T)$ be the characteristic polynomial of
Frobenius for the Jacobian
of the curve associated with $k_{0}$.  Suppose that the
distinct roots of $P(T)$ modulo $\ell$ (possibly in some extension of
$\Z/\ell\Z$) consist of $\lambda_{0}, \lambda_{1}, \ldots,
\lambda_{n}$ such that for all $i \neq j$, $\lambda_{i}\lambda_{j}
\neq 1$.  Suppose further that if any $\lambda_{i}=1$, $\lambda_{i}$
is at most a double root of $P(T)$ modulo $\ell$, and if any
$\lambda_{i}=-1$, $\lambda_{i}$ is only a simple root of $P(T)$ modulo
$\ell$.  Then there are no unramified infinite powerful pro-$\ell$
extensions of $k_{n}$, Galois over $k_0$, with no constant field
extension.
\end{theorem}

In the paper~\cite{AH}, Jeffrey Achter and the author address the
question of how many function fields are associated with a given
$P(T)$ modulo $\ell$, and thus how many fall under the purview of
Theorem~\ref{thm2} and Theorem~\ref{symm-free}.  In this paper we will
address the different but related question of how many isogeny classes
of abelian varieties have a given characteristic polynomial $P(T)$
modulo $\ell$.  As an application we find the proportion of isogeny
classes of abelian varieties with a rational point of order $\ell$.

We have chosen the following way to address these
questions, starting with the application to rational points.  Fix
distinct primes $p$ and $\ell$.  For each $r$, let $\F_{p^{r}}$ be the
finite field with $p^{r}$ elements.  By the work of Tate and Honda,
two abelian varieties are isogenous if and only if they have the same
zeta function.  Thus to each isogeny class of abelian varieties
defined over $\F_{p^{r}}$ we associate the unique polynomial $P(T)$
(the \emph{Weil polynomial} or \emph{Weil $q$-polynomial}) which is the
characteristic polynomial for the action of Frobenius and the
reciprocal of the numerator of the zeta function of any variety in the
isogeny class.  Then $\ell$ does not divide $P(1)$ if and only if each
abelian variety in the class has an $\F_{p^{r}}$-rational point of
order $\ell$.  For each $g$, there are finitely many isogeny classes
of abelian varieties with dimension $g$.  Let $d_{r,g}$ be the
fraction of isogeny classes of dimension $g$ over $\F_{p^{r}}$ for
which $\ell$ does not divide $P(1)$.  Then

\begin{theorem} \label{thm2lim}
For fixed $g$,
$$\lim_{r \to \infty} d_{r,g} =
\frac{\ell-1}{\ell}.$$
\end{theorem}

This result and the other major result of the paper could also be
obtained using the techniques of~\cite{AH}.  The proofs given here
are perhaps more elementary, and also give some access to the number
of isogeny classes and not merely the proportion satisfying each
condition.

\section{Lattices}

The proof of Theorem~\ref{thm2lim} relies on the method of counting
abelian varieties introduced by DiPippo and Howe in~\cite{DH}.  Let
$q=p^{r}$ and $I(q,g)$ be the number of isogeny classes of
$g$-dimensional abelian varieties over $\F_{q}$.  Let $P(T)$ be as
before.  Write
    $$P(T) = \prod_{j=1}^{2g}(T - \alpha_{j}).$$
Then $P(T)$ has the property that $\left|\alpha_{j}\right| = q^{1/2}$, and
the real roots, if any, have even multiplicity.  If we write
\begin{equation*}
P(T)=(T^{2g}+q^g)+a_1(T^{2g-1}+q^{g-1}T)+ \cdots + a_{g-1}(T^{g+1}+q
T^{g-1})+a_g T^g
\end{equation*}
and let $Q(T)=P(\sqrt{q}T)/q^g$, then $P(T)$ is associated with another
polynomial
$$Q(T)=(T^{2g}+1)+b_1(T^{2g-1}+T)+ \cdots + b_{g-1}(T^{g+1}+ T^{g-1})+b_g
T^g.$$

Let $V_g$ be the set of vectors $\mathbf{b}=(b_1, \ldots, b_g)$ in
$\R^{g}$ such that all of the complex roots of $Q(T)$ lie on the unit
circle and all real roots occur with even multiplicity.  Let
$\mathbf{e}_{1}, \ldots, \mathbf{e}_{g}$ be the standard basis vectors
of $\R^{g}$ and let $\Lambda_{q}$ be the lattice generated by the
vectors $q^{-i/2}\mathbf{e}_{i}$.  DiPippo and Howe explain that if
$P(T)$ is the Weil polynomial of an isogeny class then the
coefficients $a_i$ are such that $(a_{1}q^{-1/2}, \ldots,
a_{g}q^{-g/2}) \in \Lambda_{q} \cap V_{g}$.  Further, if we let
$\Lambda'_{q}$ be the lattice generated by the vectors
$q^{-1/2}\mathbf{e}_{1},\ldots, q^{-(g-1)/2}\mathbf{e}_{g-1}$ and
$pq^{-g/2}\mathbf{e}_{g}$, then all of the polynomials $P(T)$ with
coefficients $a_i$ are such that $(a_{1}q^{-1/2}, \ldots,
a_{g}q^{-g/2}) \in ( \Lambda_{q} \setminus \Lambda'_{q}) \cap V_{g}$
are exactly the Weil polynomials of isogeny classes of ordinary
varieties.  Finally, let $\Lambda''_{q}$ be the lattice generated by
the vectors $q^{-1/2}\mathbf{e}_{1},\ldots,
q^{-(g-1)/2}\mathbf{e}_{g-1}$ and $sq^{-g/2}\mathbf{e}_{g}$, where $s$
is the smallest power of $p$ such that $q$ divides $s^{2}$.  Then the set
of polynomials $P(T)$ with coefficients $a_i$ such that
$(a_{1}q^{-1/2}, \ldots, a_{g}q^{-g/2}) \in \Lambda''_{q} \cap V_{g}$
contains (perhaps properly) the set of Weil polynomials of isogeny
classes of non-ordinary varieties.

These facts are relevant because of Proposition~2.3.1 of~\cite{DH}.
In a slightly generalized form, the proposition says:

\begin{prop}[see 2.3.1 of~\cite{DH}] Let $n>0$ be an integer and let
$\Lambda \subseteq \R^{n}$ be a rectilinear lattice (possibly shifted)
with mesh $d$ at most $D$.  Then we have $$\left| \#(\Lambda \cap V_{g})
- \frac{\text{volume~} V_{n}}{\text{covolume~} \Lambda} \right| \leq
c(n,D)
\frac{d}{\text{covolume~}\Lambda}$$
for some  constant $c(n,D)$ depending only on $n$ and $D$ which can
be explicitly computed.  (We will not need the explicit computation in
this paper.)
\end{prop}

Let $v_{n}$ be the volume of $V_{n}$; Proposition~2.2.1 of~\cite{DH}
calculates it explicitly but we will not need that here.  Let
$r(q)=1-1/p$.  The lattice $\Lambda_{q}$ has covolume $q^{-g(g+1)/4}$
and mesh $q^{-1/2}$.  The lattice $\Lambda'_{q}$ has covolume
$pq^{-g(g+1)/4}$, and it has mesh $q^{-1/2}$ unless $g=2$ and $q=p$,
in which case it has mesh $1$.  Lastly, the lattice $\Lambda''_{q}$
has covolume $sq^{-g(g+1)/4}$ and its mesh is at most $1$.  It is then
an easy consequence of the proposition that
\begin{multline*}
v_{g}r(q)q^{g(g+1)/4} -
2 c(g,1)  q^{g(g+1)/4-1/2} \\
\begin{split}
&\leq
I(q,g)    \\
&\leq
v_{g} r(q) q^{g(g+1)/4} + \left(v_{g}+
3  c(g,1) \right)q^{g(g+1)/4-1/2}.
\end{split}
\end{multline*}
(See~\cite{DH} for details.)

Now let $I_{\ell}(q,g)$ be the number of isogeny classes of
$g$-dimensional abelian varieties over $\F_{q}$ such that $\ell$
divides $P(1)$.  Using the above notation we have
$$P(1)=(1+q^g)+a_1(1+q^{g-1})+ \cdots + a_{g-1}(1+q)+a_g.$$
Then
$$I_{\ell}(q,g) = \sum_{\substack{(1+q^g)+m_1(1+q^{g-1})+ \cdots +
m_{g-1}(1+q)+m_g \equiv 0 \pmod{\ell} \\
0\leq m_{i} < \ell}} I_{m_{1}, \ldots, m_{g}}(q,g)$$
where $I_{m_{1}, \ldots, m_{g}}(q,g)$ is the number of isogeny classes of
$g$-dimensional abelian varieties over $\F_{q}$ such that $a_{i}
\equiv m_{i}$ modulo $\ell$.  There are exactly $\ell^{g-1}$ terms on
the right hand side of this expression.

Now let $\Lambda_{m_{1}, \ldots, m_{g}}$ be the lattice generated by
the vectors $\ell q^{-i/2}\mathbf{e}_{i}$ and then shifted by
$\sum_{i} m_{i} q^{-i/2}\mathbf{e}_{i}$, and let $\Lambda'_{m_{1},
\ldots, m_{g}} = \Lambda_{m_{1}, \ldots, m_{g}} \cap \Lambda'_{q}$ and
$\Lambda''_{m_{1}, \ldots, m_{g}} = \Lambda_{m_{1}, \ldots, m_{g}}
\cap \Lambda''_{q}$.  Then $\Lambda_{m_{1}, \ldots, m_{g}}$ has
covolume $\ell^{g} q^{-g(g+1)/4}$ and mesh $\ell q^{-1/2}$;
$\Lambda'_{m_{1}, \ldots, m_{g}}$ has covolume
$\ell^{g}pq^{-g(g+1)/4}$, and it has mesh $\ell q^{-1/2}$ unless $g=2$
and $q=p$, in which case it has mesh $\ell$; and $\Lambda''_{m_{1},
\ldots, m_{g}}$ has covolume $\ell^{g}sq^{-g(g+1)/4}$ and mesh at most
$\ell$.

We can then prove:

\begin{prop} \label{ImProp}
\begin{multline*}
v_{g}r(q)q^{g(g+1)/4}\ell^{-g} -
2  c(g,\ell)q^{g(g+1)/4-1/2} \ell^{1-g}\\
\begin{split}
&\leq
I_{m_{1}, \ldots, m_{g}}(q,g)    \\
&\leq
v_{g} r(q) q^{g(g+1)/4}\ell^{-g} +
\left(v_{g}+
3 c(g,\ell)  \right)
q^{g(g+1)/4-1/2} \ell^{1-g},
\end{split}
\end{multline*}
\end{prop}

and thus:

\begin{prop}
\begin{multline*}
v_{g}r(q)q^{g(g+1)/4}\ell^{-1} -
2 c(g,\ell)  q^{g(g+1)/4-1/2} \\
\begin{split}
&\leq
I_{\ell}(q,g)    \\
&\leq
v_{g} r(q) q^{g(g+1)/4}\ell^{-1} +
\left(v_{g}+ 3 c(g,\ell) \right)
q^{g(g+1)/4-1/2}.
\end{split}
\end{multline*}
\end{prop}

Combining this with our earlier result, we get
\begin{multline*}
\frac{v_{g}r(q)q^{g(g+1)/4}\ell^{-1} - 2
c(g,\ell) q^{g(g+1)/4-1/2}}{v_{g}
r(q) q^{g(g+1)/4} +
\left(v_{g}+ 3 c(g,1) \right)
q^{g(g+1)/4-1/2}} \\
\begin{split}
&\leq
\frac{I_{\ell}(q,g)}{I(q,g)}    \\
&\leq \frac{v_{g} r(q) q^{g(g+1)/4}\ell^{-1} +
\left(v_{g}+ 3 c(g,\ell)  \right)
q^{g(g+1)/4-1/2}}{v_{g}
r(q)q^{g(g+1)/4} -
2 c(g,1)  q^{g(g+1)/4-1/2}}.
\end{split}
\end{multline*}
Thus we have:

\begin{theorem}
For fixed $g$,
$$\lim_{r \to \infty} \frac{I_{\ell}(p^{r}, g)}{I(p^{r}, g)} =
\frac{1}{\ell}.$$
\end{theorem}

from which Theorem~\ref{thm2lim} follows immediately.

\section{The general case}

Obviously, an identical argument could be used to establish the
fraction of isogeny classes of dimension $g$ for which $P(x) \equiv y$
modulo $\ell$ for any $x$ and $y$ in $\Z$.  More generally, we can establish
the fraction of isogeny classes of dimension $g$ for which $P(T)
\equiv f(T)$ modulo $\ell$ for any given polynomial $f(T)$ of the
correct form.  Fix
$$f(T) = (T^{2g}+q^g)+m_1(T^{2g-1}+q^{g-1}T)+ \cdots + m_{g-1}(T^{g+1}+q
T^{g-1})+m_g T^g.$$
For
fixed $p$ and $\ell$, let $e_{r,g}$ be the fraction of isogeny classes of
$g$-dimensional abelian varieties over $\F_{p^r}$ such that $P(T)
\equiv f(T)$ modulo $\ell$.

\begin{theorem} \label{generallim}
For fixed $g$,
$$\lim_{r \to \infty} e_{r,g} =
\frac{1}{\ell^{g}}.$$
\end{theorem}

\begin{proof}

We can follow the same argument as we did for Theorem~\ref{thm2lim}.
Let $J_\ell(q,g) =
e_{r,g} I(q,g)$ be the number of isogeny classes of
$g$-dimensional abelian varieties over $\F_{p^r}=\F_q$ such that $P(T)
\equiv f(T)$ modulo $\ell$.
Then our bounds on $J_\ell(q,g) = I_{m_{1}, \ldots, m_{g}}(q,g)$ and
$I(q,g)$ give us
\begin{multline*}
\frac{v_{g}r(q)q^{g(g+1)/4}\ell^{-g}  -
2 c(g,\ell) q^{g(g+1)/4-1/2} \ell^{1-g}}%
{v_{g} r(q)
q^{g(g+1)/4} +
\left(v_{g}+ 3 c(g,1) \right)
q^{g(g+1)/4-1/2}}\\
\begin{split}
&\leq
\frac{J_\ell(q,g)}{I(q,g)}    \\
&\leq
\frac{v_{g} r(q) q^{g(g+1)/4}\ell^{-g}  +
\left(v_{g}+ 3 c(g,\ell) \right)
q^{g(g+1)/4-1/2} \ell^{1-g}}%
{v_{g}
r(q)q^{g(g+1)/4} -
2 c(g,1) q^{g(g+1)/4-1/2}}.
\end{split}
\end{multline*}

On taking the limit, the theorem follows.
\end{proof}

\newcommand{\SortNoop}[1]{}

\end{document}